\documentclass[12pt,reqno]{amsart}

\usepackage[arrow,matrix,curve]{xy}

\usepackage[dvips]{graphicx} 

\usepackage{amssymb, latexsym, amsmath, amscd, array,
}

\theoremstyle{definition}

\numberwithin{equation}{section}

\newcommand\R {{\mathbb R}}

\newcommand\st{{\rm st}} 

\newcommand\RRR{\mbox{I\!I\!R}}

\newcommand\Los{{\L}o{\'s}}

\newcommand\Horvath{{Horv\'ath}}

\DeclareMathOperator{\adequal}{\,{}_{\ulcorner\!\urcorner}\,}

\begin{document}


\thispagestyle{empty}

\title{Leibniz's laws of continuity and homogeneity}

\author{Mikhail G. Katz
%
%
and David M. Sherry}

\address{Department of Mathematics, Bar Ilan University, Ramat Gan
52900 Israel} \email{katzmik@macs.biu.ac.il}

\address{Department of Philosophy, Northern Arizona University,
Flagstaff, AZ 86011, USA} \email{david.sherry@nau.edu}


\subjclass[2000]{Primary 
26E35;       
Secondary 
01A85,       
03A05       
}

\keywords{continuum; Hewitt; infinitesimal; law of continuity; law of
homogeneity; Leibniz; \Los; Robinson; standard part; Stevin; transfer
principle}

\bigskip\bigskip\bigskip\noindent
\begin{abstract}
We explore Leibniz's understanding of the differential calculus, and
argue that his methods were more coherent than is generally
recognized.  The foundations of the historical infinitesimal calculus
of Newton and Leibniz have been a target of numerous criticisms.  Some
of the critics believed to have found logical fallacies in its
foundations.  We present a detailed textual analysis of Leibniz's
seminal text \emph{Cum Prodiisset}, and argue that Leibniz's system
for differential calculus was free of contradictions.
\end{abstract}

\maketitle

\tableofcontents

\section{From the {\em characteristica universalis\/} to ideal entities}

Leibniz envisioned the creation of a universal language, {\em
characteristica universalis\/}, ambitiously designed to serve as the
vehicle for deduction and discovery in all fields of knowledge.
Couturat \cite{Co01} pointed out in 1901 that in Leibniz's vision, the
infinitesimal calculus was but the first salvo, or sample, of his {\em
characteristica universalis\/}.

Leibniz's vision of the ideal nature of mathematical entities was
remarkably modern.  His description of infinitesimals as fictional
entities shocked his disciples J. Bernoulli, l'H\^opital, and
Varignon.  And his infinitesimals certainly appeared as ``Mysteries''
to critics such as Berkeley, whose empiricist philosophy tolerated no
conceptual innovations, like infinitesimals, without an empirical
counterpart or referent:
\begin{quote}
Yet some there are, who, though they shrink at all other Mysteries,
make no difficulty of their own, who strain at a Gnat and swallow a
Camel (Berkeley 1734, \cite[section~XXXIV]{Be}).
\end{quote}

While today we are puzzled by Berkeley's rigid rejection of the idea
of an infinitely divisible continuum, he also articulated a
specifically \emph{logical} criticism of the calculus (see Sherry
1987~\cite{She87}), alleging that the system suffered from logical
flaws and even contradictions.  Even today, many historians believe
Berkeley's criticism to have been on target.  Not even Robinson
escaped this trend, praising Berkeley's criticism of the foundations
of the calculus as ``a brilliant exposure of their logical
inconsistencies'' (Robinson 1966 \cite[p.~280]{Ro66}).

We argue that, contrary to Berkeley's view, Leibniz's system for the
differential calculus was robust and free of contradiction.  Leibniz
articulated a set of coherent heuristic procedures for his calculus.
Thus, Leibniz's system incorporated versatile heuristic principles
such as his law of continuity and laws of homogeneity, which were
amenable, in the ripeness of time, to implementation as general
principles governing the manipulation of modern infinitesimal and
infinitely large quantities, such as the transfer principle and the
standard part principle.  Kanovei~\cite{Kan} and others performed
similar reconstructions of Euler's work.

We will draw on Leibniz's work, more specifically his \emph{Cum
Prodiisset}, so as to argue for the consistency of Leibniz's system
for the differential calculus.%
\footnote{The text \emph{Cum Prodiisset} sheds more light on
foundational issues than the terser 1684 text \emph{Nova methodus pro
maximis et minimis \ldots} \cite{Le84}, so named because the
determination of maxima and minima was one of the central problems of
analysis as practiced at the time, with leading scholars having
composed works of similar titles (see e.g., Fermat \cite{Fer}).}
We will also draw on the work of Leibniz historians Bos, Ferraro,
Horv\'ath, Knobloch, and Laugwitz.

Berkeley's attack on infinitesimal calculus focused specifically on
the product rule as well as the derivation of polynomials.  Once the
logical contradiction, alleged by Berkeley, is resolved in the context
of the product rule (see Section~\ref{prod}), it can similarly be
resolved (namely by applying the transcendental law of homogeneity) in
all other contexts.

In a seminal 1974 study of Leibnizian methodology, Bos described a
\emph{pair} of distinct approaches to justifying the calculus:
\begin{quote}
Leibniz considered two different approaches to the foundations of the
calculus; one connected with the classical methods of proof by
``exhaustion", the other in connection with a law of continuity (Bos
\cite[item 4.2, p.~55]{Bos}).
\end{quote}
The first approach relies on an Archimedean ``exhaustion''
methodology.  We will therefore refer to it as the A-methodology.  The
other methodology exploits infinitesimals and the law of continuity.
We will refer to it as the B-methodology, in an allusion to Johann
Bernoulli, who, having learned an infinitesimal methodology from
Leibniz, never wavered from it.

Knobloch explains the role of Leibniz's law of continuity in the
following terms:
\begin{quote}
In his treatise Leibniz used a dozen rules which constitute his
arithmetic of the infinite.  He just applied them without
demonstrating them, only relying on the {\em law of continuity\/}: The
rules of the finite remains valid in the domain of the infinite
(Knobloch \cite[p.~67]{Kn}).
\end{quote}
Laugwitz pointed out that Leibniz's law of continuity
\begin{quote}
contains an {\em a priori\/} assumption: our mathematical universe of
discourse contains both finite objects and infinite ones (Laugwitz
1992, \cite[p.~145]{Lau92}).
\end{quote}
What is the ontological status of such infinitary (infinitesimal or
infinite) objects in Leibniz's theory?  Leibniz's was a remarkably
modern insight that mathematical entities need not have a {\em
referent\/}, or empirical counterpart.  The fictional nature of
infinitesimals was stressed by Leibniz in 1706:
\begin{quote}
Philosophically speaking, I no more admit magnitudes infinitely small
than infinitely great \ldots I take both for mental fictions, as more
convenient ways of speaking, and adapted to calculation, just like
imaginary roots are in algebra. 
%
(Leibniz to Des Bosses, 11 march 1706; in Gerhardt \cite[II,
p.~305]{GP})
\end{quote}

Infinitesimals, like imaginaries, were well-founded fictions to
Leibniz.  The nature of Leibniz's infinitesimals is further clarified
by Ferraro, who analyzes a lengthy quotation from Leibniz's famous
letter \cite{Le02} to Varignon of 1702 (which we do not reproduce here
to save space), in the following terms:
\begin{quote}
According to Leibniz, imaginary numbers, infinite numbers,
infinitesimals, the powers whose exponents were not ``ordinary"
numbers and other mathematical notions are not mere inventions; they
are auxiliary and ideal quantities that [\ldots] serve to shorten the
path of thought (Ferraro \cite[p.~35]{Fe08}).
\end{quote}
On Ferraro's view, Leibniz's infinitesimals enjoy an {\em ideal\/}
ontological status similar to that of the complex numbers, {\em
surd\/} (irrational) exponents, and other ideal quantities.  

In the next section, we will examine Leibniz's foundational stance as
expressed in his seminal text \emph{Cum Prodiisset}.

\section{\textbf{\emph{Cum Prodiisset}}}
\label{CP}

Leibniz's text \emph{Cum Prodiisset} \cite{Le01c} (translated by Child
\cite{Ch}) dates from around 1701 according to modern scholars.  The
text is of crucial importance in understanding Leibniz's foundational
stance. We will analyze it in detail in this section.

\subsection{Law of Continuity, with examples}
\label{42}
\label{4.2}

Leibniz formulates his law of continuity in the following terms:
\begin{quote}
Proposito quocunque transitu continuo in aliquem terminum desinente,
liceat raciocinationem communem instituere, qua ultimus terminus
comprehendatur (Leibniz \cite[p.~40]{Le01c}).
\end{quote}
The passage can be translated as follows:
\begin{quote}
In any supposed continuous transition, ending in any terminus, it is
permissible to institute a general reasoning, in which the final
terminus may also be included.
\end{quote}
We have deliberately avoided using the term \emph{limit} in our
translation.%
\footnote{This is consonant with Child's translation: ``In any
supposed transition, ending in any terminus, it is permissible to
institute a general reasoning, in which the final terminus may also be
included" \cite[p.~147]{Ch}.  We have reinstated the adjective
\emph{continuous} modifying \emph{transition} (deleted by Child
possibly in an attempt to downplay a perceived logical circularity of
defining a law of continuity in terms of continuity itself).}

In fact, translating \emph{terminus} as \emph{limit} would involve a
methodological error, because the term \emph{limit} misleadingly
suggests its modern technical meaning of a real-valued operation
applied to sequences or functions.  In a similar vein, Bos notes that

\begin{quote}
the fundamental concepts of the Leibnizian infinitesimal calculus can
best be understood as extrapolations to the actually infinite of
concepts of the calculus of finite sequences.  I use the term
``extrapolation" here to preclude any idea of taking a limit (Bos
\cite[p.~13]{Bos}).%
\footnote{Bos goes on specifically to criticize the Bourbaki's
\emph{limite} wording ``(Leibniz) se tient tr\`es pr\`es du calcul des
diff\'erences, dont son calcul diff\'erentiel se d\'eduit par un
passage \`a la limite" (Bourbaki \cite[p.~208]{Bou}).}
\end{quote}

Leibniz gives several examples of the application of his Law of
Continuity.  We will focus on the following three examples.  

\begin{enumerate}
\item
In the context of a discussion of parallel lines, he writes:
\begin{quote}
when the straight line BP ultimately becomes parallel to the straight
line VA, even then it converges toward it or makes an angle with it,
only that the angle is then infinitely small (Child
\cite[p. 148]{Ch}).
\end{quote}
\item
Invoking the idea that the term equality may refer to equality up to
an infinitesimal error, Leibniz writes: 
\begin{quote}
when one straight line%
\footnote{Here Leibniz is using the term \emph{line} in its generic
meaning of a \emph{segment}.}
is equal to another, it is said to be unequal to it, but that the
difference is infinitely small \cite[p. 148]{Ch}.%
\footnote{Equality up to an infinitesimal is a state of transition
from inequality to equality (this anticipates the transcendental law
of homogeneity dealt with in Section~\ref{homogeneity}).}
\end{quote}
\item
Finally, a conception of a parabola expressed by means of an ellipse
with an infinitely removed focal point is evoked in the following
terms:
\begin{quote}
a parabola is the ultimate form of an ellipse, in which the second
focus is at an infinite distance from the given focus nearest to the
given vertex \cite[p. 148]{Ch}.
\end{quote}
\end{enumerate}

\subsection{\emph{Status transitus}}
\label{status}

We return to our analysis of the law of continuity as formulated in
\emph{Cum Prodiisset}.  Leibniz introduces his next observation by the
clause ``of course it is really true that'', and notes that ``straight
lines which are parallel never meet" \cite[p. 148]{Ch}; that ``things
which are absolutely equal have a difference which is absolutely
nothing" \cite[p. 148]{Ch}; and that ``a parabola is not an ellipse at
all" \cite[p .149]{Ch}.  These remarks would seem to rule out the
wondrous entities of the previous subsection.  How does one, then,
account for these examples?  Leibniz provides an explanation in terms
of a state of transition (\emph{status transitus} in the original
Latin \cite[p. 42]{Le01c}):
\begin{quote}
a state of transition may be imagined, or one of evanescence, in which
indeed there has not yet arisen exact equality \ldots{} or
parallelism, but in which it is passing into such a state, that the
difference is less than any assignable quantity; also that in this
state there will still remain some difference, \ldots{} some angle,
but in each case one that is infinitely small; and the distance of the
point of intersection, or the variable focus, from the fixed focus
will be infinitely great, and the parabola may be included under the
heading of an ellipse \cite[p.~149]{Ch}.
\end{quote}

A state of transition in which ``there has not yet arisen exact
equality" refers to example (2) in Subsection~\ref{4.2};
``parallelism" refers to example (1); including parabola under the
heading of ellipse is example (3).  

Thus, the term \emph{terminus} encompasses the \emph{status
transitus}, involving a passage into an assignable entity, while being
as yet inassignable.  Translating \emph{terminus} as limit amounts to
translating it as an assignable entity, the antonym of the meaning
intended by Leibniz.

The observation that Leibniz's \emph{status transitus} is an
inassignable quantity is confirmed by Leibniz's conceding that its
metaphysical status is ``open to question'':

\begin{quote}
whether such a state of instantaneous transition from inequality to
equality, \ldots{} from convergence [i.e., lines meeting--the authors]
to parallelism, or anything of the sort, can be sustained in a
rigorous or metaphysical sense, or whether infinite extensions
successively greater and greater, or infinitely small ones
successively less and less, are legitimate considerations, is a matter
that I own to be possibly open to question \cite[p. 149]{Ch}.
\end{quote}

Yet Leibniz asserts that infinitesimals may be utilized independently
of metaphysical controversies:
\begin{quote}
but for him who would discuss these matters, it is not
necessary to fall back upon metaphysical controversies,
such as the composition of the continuum, or to make
geometrical matters depend thereon \cite[p. 149-150]{Ch}.
\end{quote}
To summarize, Leibniz holds that the inassignable status of
\emph{status transitus} is no obstacle to its effective use in
geometry.  The point is reiterated in the next paragraph:

\begin{quote}
If any one wishes to understand these [i.e. the infinitely great or
the infinitely small--the authors] as the ultimate things, or as truly
infinite, it can be done, and that too without falling back upon a
controversy about the reality of extensions, or of infinite continuums
in general, or of the infinitely small, ay, even though he think that
such things are utterly impossible; it will be sufficient simply to
make use of them as a tool that has advantages for the purpose of the
calculation, just as the algebraists retain imaginary roots with great
profit \cite[p. 150]{Ch}.
\end{quote}
Leibniz has just asserted the possibility of the \emph{mathematical}
infinite: ``it can be done", without \emph{philosophical} commitments
as to its ontological reality.

\subsection{Mathematical implementation of \emph{status transitus}}
\label{45}

We will illustrate Leibniz's concept of \emph{status transitus} by
implementing it mathematically in the three examples mentioned by
Leibniz.  

In Subsection~\ref{status}, we mentioned that Leibniz viewed
infinitesimals as fictions, and so his methods avoided any
metaphysical commitments.  Ishiguro (1990, \cite{Is}) and others
viewed a Leibnizian infinitesimal as a \emph{logical} fiction,
involving a syncategorematic paraphrase with a hidden quantifier
applied to ordinary real values.  We have argued against the
syncategorematic interpretation in~\cite{KS}.  Rather, Leibnizian
infinitesimals are \emph{pure} fictions akin to imaginaries.  Thus,
Leibniz exploited infinitely large and infinitely small quantities in
the same sense in which Albert Girard (1595-1632) and others exploited
imaginary roots in order to simplify algebra.  In neither case is
there a commitment to corresponding mathematical entities.  Thus,
Leibniz anticipates modern formalist positions such as Hilbert's and
Robinson's.

Of course, the structural properties of Leibniz's infinite and
infinitely small quantities are different from those of modern day
infinitesimals.  Nonetheless, modern theories of infinitesimals
\emph{are} a way of implementing Leibniz's heuristic procedures.
Thus, Example (2) can be illustrated as follows.  Leibniz denotes a
finite positive quantity by
\[
(d)x
\]
(Bos \cite[p.~57]{Bos} replaced this by~$\underline{\text{d}}x$).  The
assignable quantity~$(d)x$ passes via infinitesimal~$dx$ on its way to
absolute~$0$.  Then the infinitesimal~$dx$ is the \emph{status
transitus}.  Zero is merely the assignable \emph{shadow}%
\footnote{See footnote~\ref{sha}.}
of the infinitesimal.  Then a \emph{line} (i.e., segment) of
length~$2x+dx$ will be equal to one of length~$2x$, up to an
infinitesimal.  This particular \emph{status transitus} is the
foundation rock of the Leibnizian definition of the differential
quotient.

Example (1) of parallel lines can be elaborated as follows.  Let us
follow Leibniz in building the line through~$(0, 1)$ parallel to the
$x$-axis in the plane.  Line~$L_H$ with~$y$-intercept~$1$ and
$x$-intercept~$H$ is given by~$y = 1 -\frac{x}{H}$.  For infinite~$H$,
the line~$L_H$ has negative infinitesimal slope, meets the~$x$-axis at
an infinite point, and forms an infinitesimal angle with the~$x$-axis
at the point where they meet.  We will denote by~$\st(x)$ the
assignable (i.e., real) shadow%
\footnote{\label{sha}Here ``st" stands for the standard part function
in the context of the hyperreals.  Of course, Leibniz used neither the
term ``shadow'' nor the ``st'' notation.  Rather, these notations from
modern infinitesimal analysis implement mathematically a heuristic
principle of Leibniz's called the transcendental law of homogeneity,
discussed in Section~\ref{homogeneity}.}
 of a finite~$x$.  Then every finite point~$(x,y)\in L_H$ satisfies
\[
\begin{aligned}
\st(x, y) & = (\st(x), \st(y)) \\& = \left(\st(x), \st\left(1
-\tfrac{x}{H}\right)\right) \\& = (\st(x), 1).
\end{aligned}
\]
Hence the finite portion of~$L_H$ is infinitely close to the
line~$y=1$ parallel to the~$x$-axis, which is its \emph{shadow}. Thus,
the parallel line is constructed by varying the oblique line depending
on a parameter. Such variation passes via the \emph{status transitus}
defined by an infinite value of~$H$.

To implement example (3), let's follow Leibniz in deforming an
ellipse, via a \emph{status transitus}, into a parabola. The ellipse
with vertex (apex) at~$(0,-1)$ and with foci at the origin and
at~$(0;H)$ is given by
\begin{equation}
\label{4.1}
\sqrt{x^2 + y^2} +\sqrt{x^2 + (y - H)^2} = H + 2
\end{equation}
We square \eqref{4.1} to obtain
\begin{equation}
\label{4.2b}
x^2 +y^2 +x^2 +(H-y)^2 + 2\sqrt{(x^2+y^2)(x^2+(H-y)^2)} =H^2+4H+4
\end{equation}
We move the radical to one side
\begin{equation}
\label{4.3}
2\sqrt{(x^2+y^2)(x^2+(H-y)^2)} =H^2+4H+4
-\left(x^2+y^2+x^2+(H-y)^2\right)
\end{equation}
and square again.  After cancellation we see that \eqref{4.1} is
equivalent to
\begin{equation}
\label{4.4}
\left(y + 2 + \tfrac{2}{H}\right)^2 - (x^2 + y^2)\left(1 +
\tfrac{4}{H} + \tfrac{4}{H^2}\right) = 0.
\end{equation}
The calculation \eqref{4.1} through \eqref{4.4} depends on habits of
\emph{general reasoning} such as:
\begin{itemize}
\item
squaring undoes a radical;
\item
the binomial formula;
\item
terms in an equation can be transfered to the other side; etc.
\end{itemize}
\emph{General reasoning} of this type is familiar in the realm of
ordinary assignable (finite) numbers, but why does it remain valid
when applied to the, fictional, ``realm'' of inassignable (infinite or
infinitesimal) numbers?  The validity of transfering such
\emph{general reasoning} originally \emph{instituted} in the finite
realm, to the ``realm'' of the infinite is precisely the content of
Leibniz's law of continuity.%
\footnote{When the \emph{general reasoning} being transfered to the
infinite ``realm'' is generalized to encompass arbitrary elementary
properties (i.e. first order properties), one obtains the
\Los-Robinson transfer principle.}

We therefore apply Leibniz's law of continuity to equation~\eqref{4.4}
for an infinite~$H$.  The resulting entity is still an ellipse of
sorts, to the extent that it satisfies all of the equations
\eqref{4.1} through \eqref{4.4}. However, this entity is no longer
finite. It represents a Leibnizian \emph{status transitus} between
ellipse and parabola. This \emph{status transitus} has foci at the
origin and at an infinitely distant point~$(0,H)$.  Assuming~$x$ and
$y$ are finite, we set~$x_0 = \st(x)$ and~$y_0 = \st(y)$, to obtain an
equation for a real shadow of this entity:
\[
\begin{aligned}
\st & \left( \left(y + 2 + \tfrac{2}{H} \right)^2 - (x^2 + y^2)\left(1
+ \tfrac{4}{H} + \tfrac{4}{H^2} \right) \right) = \\& = 
\left(y_0 + 2 + \st \left(\tfrac{2}{H}\right) \right)^2 - \left(x^2_0 +
y^2_0\right) \left(1 + \st \left( \tfrac{4}{H} + \tfrac{4}{H^2}
\right)\right) \\& = (y_0 + 2)^2 - \left(x^2_0 + y^2_0\right) \\& = 0.
\end{aligned}
\]
Simplifying, we obtain
\begin{equation}
\label{43}
y_0= \frac{x_0^2}{4}-1.
\end{equation}
Thus, the finite portion of the \emph{status transitus} \eqref{4.4} is
infinitely close to its \emph{shadow}~\eqref{43}, namely the real
parabola~$y= \frac{x^2}{4}-1$.  This is the kind of payoff Leibniz is
seeking with his law of continuity.

Some historians have been reluctant to interpret Leibniz's mathematics
in terms of modern mathematical theories.  Thus, Dauben presents a
list of authors, including Detlef Laugwitz, who ``have used
nonstandard analysis to rehabilitate or `vindicate' earlier
infinitesimalists'', and concludes:
\begin{quote}
Leibniz, Euler, and Cauchy [\ldots] had, in the views of some
commentators, ``Robinsonian'' nonstandard infinitesimals in mind from
the beginning.  The most detailed and methodically [sic] sophisticated
of such treatments to date is that provided by Imre Lakatos; in what
follows, it is his analysis of Cauchy that is emphasized (Dauben 1988,
\cite[p.~179]{Da88}).
\end{quote}
However, Lakatos's treatment was certainly not ``the most detailed and
methodically sophisticated'' one by the time Dauben's text appeared in
1988.  Thus, in 1987, Laugwitz had published a detailed scholarly
study of Cauchy in \emph{Historia Mathematica} (Laugwitz
\cite{Lau87}).  Laugwitz's text in \emph{Historia Mathematica} seems
to be the published version of his 1985 preprint \emph{Cauchy and
infinitesimals}.  Laugwitz's 1985 preprint does appear in Dauben's
bibliography (Dauben, 1988 \cite[p.~199]{Da88}), indicating that
Dauben was familiar with it.  It is odd to suggest, as Dauben seems
to, that a scholarly study published in \emph{Historia Mathematica}
would countenance a view that Leibniz and Cauchy could have had
```Robinsonian' nonstandard infinitesimals in mind from the
beginning''.  Surely Dauben has committed a strawman fallacy here.

To a historian, the claim that Leibniz's differential calculus was
free of logical fallacies may seem analogous to claiming that the
circle can be squared%
\footnote{Such was indeed the tenor of a recent referee report, see

\noindent
http://u.cs.biu.ac.il/$\sim$katzmik/straw2.html}
--but only if the historian embraces the triumviratist story of
analysis as an ineluctable march from incoherent infinitesimalism
toward the yawning heights of Weierstrassian epsilontics.

Rather, Lakatos, Laugwitz, Br\aa ting (2007, \cite{Br}) and others
have argued that infinitesimals as employed by Leibniz, Euler, and
Cauchy have found a set-theoretic implementation in the framework of
modern theories of infinitesimals.  The existence of such
implementations indicates that the historical infinitesimals were less
prone to contradiction than has been routinely maintained by
triumvirate historians.%
\footnote{\label{triumvirate}C.~Boyer refers to Cantor, Dedekind, and
Weierstrass as ``the great triumvirate'', see \cite[p.~298]{Boy}.}
The issue is dealt with in more detail by Katz \& Katz \cite{KK12},
\cite{KK11b}, \cite{KK11c}, \cite{KK11d}; B\l aszczyk et
al. \cite{BKS}; Borovik et al. \cite{BK}; Katz \& Leichtnam~\cite{KL};
and Katz, Schaps, \& Shnider \cite{KSS12}.

\subsection{Assignable \emph{versus} unassignable}
\label{assign}

In this section, we will retain the term ``unassignable'' from Child's
translation \cite{Ch} (\emph{inassignabiles} in the original Latin,
see \cite[p.~46]{Le01c}).  After introducing finite quantities~$(d)x,
(d)y, (d)z$, Leibniz notes that
\begin{quote}
the unassignables~$dx$ and~$dy$ may be substituted for them by a
method of supposition even in the case when they are evanescent (Child
\cite[p.~153]{Ch}).
\end{quote}
Leibniz proceeds to derive his multiplicative law in the case~$ay=xv$.
Simplifying the differential quotient, Leibniz obtains
\begin{equation}
\label{46}
\frac{ady}{dx}= \frac{xdv}{dx} + v +dv.
\end{equation}
At this point Leibniz proposes to transfer ``the matter, as we may, to
straight lines that never become evanescent'', obtaining%
\footnote{Child incorrectly transcribes formula~\eqref{47} from
Gerhardt, replacing the equality sign in Gerhardt by a plus sign.
Note that Leibniz himself used the sign~$\adequal$ (see McClenon
\cite[p.~371]{Mc23}).}
\begin{equation}
\label{47}
\frac{a\,(d)y}{(d)x}= \frac{x\,(d)v}{(d)x} + v +dv.
\end{equation}
The advantage of \eqref{47} over \eqref{46} is that the expressions
$\frac{(d)y}{(d)x}$ and~$\frac{(d)v}{(d)x}$ are assignable (real).  At
this stage, Leibniz points out that ``$dv$ is superfluous''.  The
reason given is that ``it alone can become evanescent''.  The
transcendental law of homogeneity (see Section~\ref{homogeneity}) is
not mentioned explicitly in \emph{Cum Prodiisset}; therefore the
discussion of this step necessarily remains a bit vague.  Discarding
the~$dv$ term, one obtains the expected product
formula~$\frac{a\,(d)y}{(d)x}= \frac{x\,(d)v}{(d)x} + v$ in this case.
Note that thinking of the left hand side of~\eqref{47} as the
assignable \emph{shadow} of the right hand side is consistent with
Leibniz's example (2) (see Subsection~\ref{42}).

\subsection{\emph{Souverain principe}}
\label{2feb02}

In a 2~feb.~1702 letter to Varignon, Leibniz formulated the law of
continuity as follows:
\begin{quote}
[\ldots] et il se trouve que les r\`egles du fini r\'eussissent dans
l'infini comme s'il y avait des atomes (c'est \`a dire des
\'el\'ements assignables de la nature) quoiqu'il n'y en ait point la
mati\`ere \'etant actuellement sousdivis\'ee sans fin; et que vice
versa les r\`egles de l'infini r\'eussissent dans le fini, comme s'il
y'avait des infiniment petits m\'etaphysiques, quoiqu'on n'en n'ait
point besoin; et que la division de la mati\`ere ne parvienne jamais
\`a des parcelles infiniment petites: c'est parce que tout se gouverne
par raison, et qu'autrement il n'aurait point de science ni r\`egle,
ce qui ne serait point conforme avec la nature du souverain principe
(Leibniz \cite[p.~350]{Le02}).
\end{quote}
This formulation was cited in (Robinson 1966 \cite[p.~262]{Ro66}).  To
summarize: \emph{the rules of the finite succeed in the infinite, and
conversely}.

\section{Assignable and inassignable quantities}
\label{homogeneity}

How did Leibniz view the relation of assignable and inassignable
quantities?  

\subsection{Relation of being infinitely close}

The rule governing infinitesimal calculation that
Knobloch represents as Leibniz's rule~2.2, states:
\begin{quote}
2.2~$x, y$ finite,~$x = (y +$ infinitely small)~$\iff$~$x - y \approx
0$ (not assignable difference) (Knobloch \cite[p.~67]{Kn}).
\end{quote}
Here the pair of parallel wavy lines represents the relation of being
infinitely close.  Leibnizian assignable quantities mark locations in
what would be called today an Archimedean continuum, or A-continuum
for short.  Such a continuum stems from the 16th century work of Simon
Stevin (1548-1620)~\cite{St85}, \cite{Ste}.  Stevin initiated a
systematic approach to decimal representation of measuring numbers,
marking a transition from a discrete arithmetic as practiced by the
Greeks, to the arithmetic of the continuum taken for granted today
(see Malet~\cite{Mal06}, Naets~\cite{Nae}, and B\l aszczyk et
al. \cite{BKS}).

Closely related to the distinction between the A- and B-methodologies
is a distinction between two types of continua, which could be called
an A-continuum and a B-continuum.  The latter encompasses inassignable
entities such as infinitesimals (in addition to assignable ones), and
can be described as a ``thick'' continuum.%
\footnote{The B-continuum can be thought of as ``thicker'' than the
A-continuum because the B-continuum is, as it were, packed chock-full
of numbers, including infinitesimals.}
On occasion, Leibniz describes such entities as ``incomparable
quantities'', and defines them in terms of the violation of what today
is called the Archimedean property.  Thus, Leibniz writes in a letter
to l'H\^opital:
\begin{quote}
I call incomparable quantities of which the one can not become larger
than the other if multiplied by any finite number.  This conception is
in accordance with the fifth definition of the fifth book of Euclid
(Leibniz \cite[p.~288]{Le95a}, cited in \Horvath~\cite[p.~63]{Ho86}).%
\footnote{\Horvath{} notes that Leibniz is actually referring to the
{\em fourth\/} definition of the fifth book.}
\end{quote}

\subsection{Transcendental law of homogeneity}

To mediate between assignable and inassignable quantities, Leibniz
developed an additional principle called the {\em transcendental law
of homogeneity\/}.  Leibniz's transcendental law of homogeneity
governs equations involving differentials.  Bos interprets it as
follows:
\begin{quote}
A quantity which is infinitely small with respect to another quantity
can be neglected if compared with that quantity.  Thus all terms in an
equation except those of the highest order of infinity, or the lowest
order of infinite smallness, can be discarded.  For instance,
\begin{equation}
\label{adeq}
a+dx =a
\end{equation}
\[
dx+ddy=dx
\]
etc.  The resulting equations satisfy this [\dots] requirement of
homogeneity (Bos \cite[p.~33]{Bos} paraphrasing Leibniz 1710,
\cite[p.~381-382]{Le10b}).
\end{quote}
The title of Leibniz's 1710 text is \emph{Symbolismus memorabilis
calculi algebraici et infinitesimalis in comparatione potentiarum et
differentiarum, et de lege homogeneorum transcendentali}.  The
inclusion of the transcendental law of homogeneity (\emph{lex
homogeneorum transcendentalis}) in the title of the text attests to
the importance Leibniz attached to this law.

How did Leibniz use the transcendental law of homogeneity in
developing the calculus?  In Section~\ref{prod}, we will illustrate an
application of the transcendental law of homogeneity to the particular
example of the derivation of the product rule.

\section{Justification of the product rule}
\label{prod}

The issue is the justification of the last step in the following
calculation:
\begin{equation}
\label{41}
\begin{aligned}
d(uv) &= (u+du)(v+dv)-uv=udv+vdu+du\,dv \\ & =udv+vdu.
\end{aligned}
\end{equation}

The last step in the calculation~\eqref{41}, namely
\[
{udv+vdu} + {du\,dv} = {udv+vdu}
\]
is an application of Leibniz's law of homogeneity.%
\footnote{Leibniz had two laws of homogeneity, one for dimension and
the other for the order of infinitesimalness.  Bos notes that they
`disappeared from later developments' \cite[p.~35]{Bos}, referring to
Euler and Lagrange.}

In his 1701 text {\em Cum Prodiisset\/} \cite[p.~46-47]{Le01c},
Leibniz presents an alternative justification of the product rule (see
Bos \cite[p.~58]{Bos}).  Here he divides by~$dx$ and argues with
differential quotients rather than differentials.  We analyzed
Leibniz's calculation in Subsection~\ref{assign}.  Adjusting Leibniz's
notation to fit with~\eqref{41}, we obtain an equivalent calculation%
\footnote{The special case treated by Leibniz is~$u(x)=x$.  This
limitation does not affect the conceptual structure of the argument.}
\[
\begin{aligned}
\frac{d(uv)}{dx} &= \frac{(u+du)(v+dv)-uv}{dx} \\&=
\frac{udv+vdu+du\,dv}{dx} \\&= \frac{udv+vdu}{dx} + \frac{du\,dv}{dx} \\&=
\frac{udv+vdu}{dx}.
\end{aligned}
\]
Under suitable conditions the term~$\frac{du\,dv}{dx}$ is
infinitesimal, and therefore the last step
\[
\frac{udv+vdu}{dx} + \frac{du\,dv}{dx} =
u\,\frac{dv}{dx}+v\,\frac{du}{dx}
\]
is legitimized as a special case of the transcendental law of
homogeneity, which interprets the equality sign in~\eqref{adeq} as the
relation of being infinitely close, i.e., an equality up to
infinitesimal error.  Note that the use of the equality sign ``$=$''
to denote a non-symmetric relation of discarding the ``inhomogeneous''
terms as in \eqref{adeq} and \eqref{41} should hardly shock the modern
reader used to the ``big-O'' notation: we write~$\sin x = O(1)$, but
we would certainly not write~$O(1)=\sin x$.  Leibniz's transcendental
law of homogeneity involved such an ``asymmetric'' relation, since it
replaced an inassignable quantity such as~$2x+dx$ by the assignable
quantity~$2x$.

\section{Was Leibniz's system for differential calculus, consistent?}

Berkeley's {\em logical criticism\/} of the calculus is that the
evanescent increment is first assumed to be non-zero to set up an
algebraic expression, and then {\em treated as zero\/} in {\em
discarding\/} the terms that contained that increment when the
increment is said to vanish.  The criticism, however, involves a
misunderstanding of Leibniz's method.  The rebuttal of Berkeley's
criticism is that the evanescent increment need {\em not\/} be
``treated as zero'', but, rather, merely {\em discarded\/} through an
application of the transcendental law of homogeneity by Leibniz, as
illustrated in the previous section in the case of the product rule.

\begin{figure}
\includegraphics[height=2in]{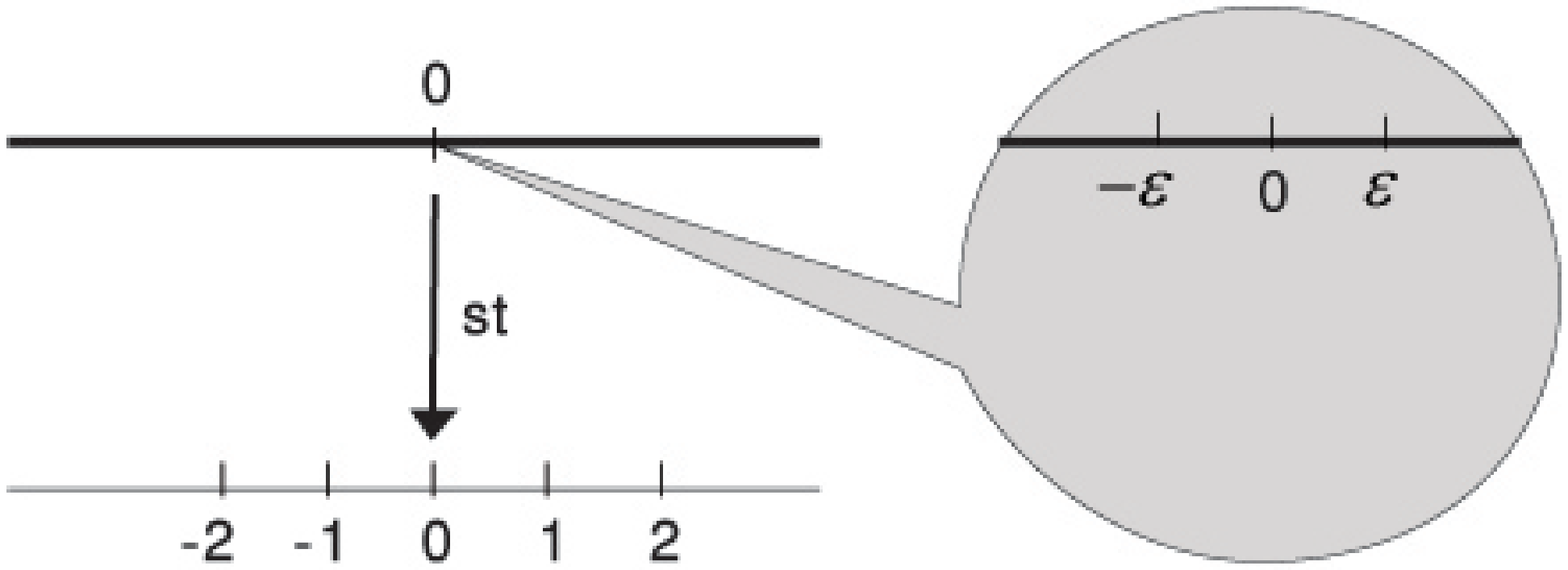}
\caption{\textsf{Zooming in on infinitesimal~$\epsilon$ (here st$(\pm
\epsilon)=0$)}}
\label{tamar}
\end{figure}

While consistent, Leibniz's system unquestionably relied on heuristic
principles such as the laws of continuity and homogeneity, and would
thus fall short of a standard of rigor \emph{if} measured by today's
criteria.  On the other hand, the consistency and resilience of
Leibniz's system is confirmed through the development of modern
implementations of Leibniz's heuristic principles.  Thus, in the
1940s, Hewitt~\cite{Hew} developed a modern implementation of a
hyperreal B-continuum extending~$\R$, by means of a technique referred
to today as the ultrapower construction.  We will denote such a
B-continuum by the new symbol~$\RRR$ (``thick-R'').  Denoting by
$\RRR_\infty$ the subset of $\RRR$ consisting of inverses of nonzero
infinitesimals, we obtain a partition $\RRR=\RRR_{<\infty} \cup
\RRR_\infty$ where $\RRR_{<\infty}$ is the complement of
$\RRR_\infty$.  We then have the standard part function 
\begin{equation}
\label{51}
\st: \RRR_{<\infty}\to\R,
\end{equation}
illustrated in Figure~\ref{tamar}.  Note that the hyperreals can be
constructed out of integers (see Borovik, Jin, \& Katz~\cite{BJK}).
The traditional quotient construction using Cauchy sequences, usually
attributed to Cantor (and is actually due to M\'eray 1869, \cite{Me}
who published three years earlier than E.~Heine), can be factored
through the hyperreals (see Giordano \& Katz \cite{GK11}).  In 1955,
\Los~\cite{Lo55} proved his celebrated theorem on ultraproducts,
implying in particular that elementary (more generally, first-order)
statements over~$\R$ are true if and only if they are true
over~$\RRR$, yielding a modern implementation of the Leibnizian law of
{\em continuity\/}.  Such a result is equivalent to what is known in
the literature as the {\em transfer principle\/} (see Keisler
\cite{Ke08}).  The map that associates to every finite element
of~$\RRR$, the real number infinitely close to it, is known in the
literature as the \emph{standard part function}~\eqref{51}
(alternatively, the {\em shadow\/}).  Such a map is a mathematical
implementation of the Leibnizian transcendental law of
\emph{homogeneity}.

\section*{Acknowledgments}

We are grateful to S. Wenmackers for helpful comments.  The influence
of Hilton Kramer (1928-2012) is obvious.


\begin{thebibliography}{AI}

\bibitem{Be} Berkeley, G.: The Analyst, a Discourse Addressed to an
Infidel Mathematician (1734).

\bibitem{BKS} B\l{}aszczyk, P.; Katz, M.; Sherry, D.: Ten
misconceptions from the history of analysis and their debunking.
\emph{Foundations of Science} (online first), see

http://dx.doi.org/10.1007/s10699-012-9285-8

and http://arxiv.org/abs/1202.4153
%
%

\bibitem{BK} Borovik, A.; Katz, M.: Who gave you the
Cauchy--Weierstrass tale?  The dual history of rigorous calculus.
\emph{Foundations of Science}, \textbf{17} (2012), no.~3, 245--276.
See http://dx.doi.org/10.1007/s10699-011-9235-x

and http://arxiv.org/abs/1108.2885

\bibitem{BJK} Borovik, A.; Jin, R.; Katz, M.: An integer construction
of infinitesimals: Toward a theory of Eudoxus hyperreals.  \emph{Notre
Dame Journal of Formal Logic} \textbf{53} (2012), no.~4, 557-570.  See
http://dx.doi.org/10.1215/00294527-1722755 and
http://arxiv.org/abs/1210.7475


\bibitem{Bos} Bos, H. J. M.: Differentials, higher-order differentials
and the derivative in the Leibnizian calculus.  {\em Arch. History
Exact Sci.\/} \textbf{14} (1974), 1--90.

\bibitem{Bou} Bourbaki, N.: El\'ements d'histoire des math\'ematiques.
\emph{Histoire de la Pens\'ee}, IV.  Hermann, Paris, 1960.

\bibitem{Boy} Boyer, C.: The concepts of the calculus.  Hafner
Publishing Company, 1949.

\bibitem{Br} Br\aa ting, K.: A new look at E. G. Bj\"orling and the
Cauchy sum theorem.  {\em Arch. Hist. Exact Sci.\/} \textbf{61}
(2007), no.~5, 519--535.
%
%

\bibitem{Ch} Child, J. M. (ed.): The early mathematical manuscripts of
Leibniz.  Translated from the Latin texts published by Carl Immanuel
Gerhardt with critical and historical notes by J. M. Child.
Chicago-London: The Open Court Publishing Co., 1920.


\bibitem{Co01} Couturat, L.: La logique de Leibniz.  F. Alcan, Paris,
1901; reprinted, Olms, Hildesheim, 1961.


\bibitem{Da88} Dauben, J.: Abraham Robinson and Nonstandard Analysis:
History, Philosophy, and Foundations of Mathematics.  In William
Aspray and Philip Kitcher, eds.  History and philosophy of modern
mathematics (Minneapolis, MN, 1985), 177--200, Minnesota
Stud. Philos. Sci., XI, Univ. Minnesota Press, Minneapolis, MN, 1988.
Available online at

http://www.mcps.umn.edu/philosophy/11\_7Dauben.pdf


\bibitem{Fer} Fermat, P.: M\'ethode pour la recherche du maximum et du
minimum. p. 121-156 in Tannery's edition \cite{Tan3}.

\bibitem{Fe08} Ferraro, G.: Geometrical quantities and series in
Leibniz.  In The Rise and Development of the Theory of Series up to
the Early 1820s Sources and Studies in the History of Mathematics and
Physical Sciences, 2008, Part I, 25-44, see
http://dx.doi.org/10.1007/978-0-387-73468-2\_2

\bibitem{Ge46} Gerhardt:, C. I. (ed.): Historia et Origo calculi
differentialis a G. G. Leibnitio conscripta, ed. C. I. Gerhardt,
Hannover, 1846.


\bibitem{Ge50} Gerhardt, C. I. (ed.): Leibnizens mathematische
Schriften (Berlin and Halle: Eidmann, 1850-1863).
%
%


\bibitem{GP} Gerhardt C. I. (ed.): G. W. Leibniz: Philosophische
Schriften.  Edited by C. I. Gerhardt. 7 vols. (1875-90).  Reprint,
Hildesheim: Georg Olms Verlag, 1962.



\bibitem{GK11} Giordano, P.; Katz, M.: Two ways of obtaining
infinitesimals by refining Cantor's completion of the reals.
Preprint, 2011,

see http://arxiv.org/abs/1109.3553


\bibitem{Hew} Hewitt, E.: Rings of real-valued continuous functions.
I.  {\em Trans. Amer. Math. Soc.\/} \textbf{64} (1948), 45--99.



\bibitem{Ho86} \Horvath, M.: On the attempts made by Leibniz to
justify his calculus.  {\em Studia Leibnitiana\/} \textbf{18} (1986),
no.~1, 60--71.
%
%


\bibitem{Is} Ishiguro, H.: Leibniz's philosophy of logic and language.
Second edition.  Cambridge University Press, Cambridge, 1990.


\bibitem{Kan} Kanovei, V.: Correctness of the Euler method of
decomposing the sine function into an infinite product. (Russian)
\emph{Uspekhi Mat. Nauk} \textbf{43} (1988), no.~4 (262), 57--81, 255;
translation in \emph{Russian Math. Surveys} \textbf{43} (1988), no.~4,
65--94.



\bibitem{KKM} Kanovei, V.; Katz, M.; Mormann, T.: Tools, Objects, and
Chimeras: Connes on the Role of Hyperreals in Mathematics.
\emph{Foundations of Science} (online first).  See
http://dx.doi.org/10.1007/s10699-012-9316-5

and http://arxiv.org/abs/1211.0244




\bibitem{KK11b} Katz, K.; Katz, M.: Cauchy's continuum.
\emph{Perspectives on Science} \textbf{19} (2011), no.~4, 426-452.
See http://arxiv.org/abs/1108.4201 and

http://www.mitpressjournals.org/doi/abs/10.1162/POSC\_a\_00047



\bibitem{KK11d} Katz, K.; Katz, M.: Meaning in classical mathematics:
Is it at odds with Intuitionism?  \emph{Intellectica} \textbf{56}
(2011), no.~2, 223--302.  See

http://arxiv.org/abs/1110.5456


\bibitem{KK12} Katz, K.; Katz, M.: A Burgessian critique of
nominalistic tendencies in contemporary mathematics and its
historiography.  \emph{Foundations of Science} \textbf{17} (2012),
no.~1, 51--89.  See http://dx.doi.org/10.1007/s10699-011-9223-1

and http://arxiv.org/abs/1104.0375




\bibitem{KK11c} Katz, K.; Katz, M.: Stevin numbers and reality,
\emph{Foundations of Science} \textbf{17} (2012), no.~2, 109-123.  See
http://arxiv.org/abs/1107.3688 and

http://dx.doi.org/10.1007/s10699-011-9228-9




\bibitem{KL} Katz, M.; Leichtnam, E.: Commuting and non-commuting
infinitesimals.  \emph{American Mathematical Monthly} (to appear).



\bibitem{KSS12} Katz, M.; Schaps, D.; Shnider, S.: Almost Equal: The
Method of Adequality from Diophantus to Fermat and Beyond.
\emph{Perspectives on Science} \textbf{21} (2013), no.~3, to appear.
See http://arxiv.org/abs/1210.7750



\bibitem{KS} Katz, M.; Sherry, D.: Leibniz's infinitesimals: Their
fictionality, their modern implementations, and their foes from
Berkeley to Russell and beyond.  \emph{Erkenntnis} (online first), see
http://dx.doi.org/10.1007/s10670-012-9370-y

and http://arxiv.org/abs/1205.0174




\bibitem{Ke08} Keisler, H. J.: The ultraproduct construction.
Proceedings of the Ultramath Conference, Pisa, Italy, 2008.


\bibitem{Kn} Knobloch, E.: Leibniz's rigorous foundation of
infinitesimal geometry by means of Riemannian sums. Foundations of the
formal sciences, 1 (Berlin, 1999).  {\em Synthese\/} \textbf{133}
(2002), no.~1-2, 59--73.
%
%

\bibitem{Lau87} Laugwitz, D.: Infinitely small quantities in Cauchy's
textbooks.  {\em Historia Math.\/} \textbf{14} (1987), no.~3,
258--274.


\bibitem{Lau92} Laugwitz, D.: Leibniz' principle and omega calculus.
[A] Le labyrinthe du continu, Colloq., Cerisy-la-Salle/Fr. 1990,
144-154 (1992).

\bibitem{Le1993} Leibniz, G. W.: (1993) De quadratura arithmetica
circuli ellipseos et hyperbolae cujus corollarium est trigonometria
sine tabulis, kritisch herausgegeben und kommentiert von Eberhard
Knobloch, G\"ottingen [Abhandlungen der Akademie der Wissenschaften in
G\"ottingen, Mathematisch-physikalische Klasse 3; 43].


\bibitem{Le2004} Leibniz, G. W.: (2004) Quadrature arithm\'etique du
cercle, de l'ellipse et de l'hyperbole, Marc Parmentier (Trans. and
Ed.)/Latin text Eberhard Knobloch (Ed.), Paris: J. Vrin.


\bibitem{Le84} Leibniz, Nova methodus pro maximis et minimis \ldots,
in {\em Acta Erud.\/}, Oct. 1684.  See Gerhardt \cite{Ge50}, V,
pp. 220-226.


\bibitem{Le95a} Leibniz to l'Hospital, 21 June, 1695, in Gerhardt
\cite{Ge50}, II, pp. 287-289.


\bibitem{Le01c} Leibniz (1701) {\em Cum Prodiisset\/}\ldots mss ``Cum
prodiisset atque increbuisset Analysis mea infinitesimalis ..." in
Gerhardt \cite{Ge46}, pp.~39--50.





\bibitem{Le02} Leibniz to Varignon, 2 Febr., 1702, in Gerhardt (see
item \cite{Ge50}) IV, pp. 91--95.

\bibitem{Le10b} Leibniz (1710) Symbolismus memorabilis calculi
algebraici et infinitesimalis in comparatione potentiarum et
differentiarum, et de lege homogeneorum transcendentali.  In Gerhardt
\cite[vol.~V, pp.~377-382]{Ge50}.

\bibitem{Lo55} {\L}o{\'s}, J.: Quelques remarques, th\'eor\`emes et
probl\`emes sur les classes d\'efi\-nissables d'alg\`ebres, in
Mathematical interpretation of formal systems, {98--113},
North-Holland Publishing Co., {Amsterdam}, {1955}.
%
%


\bibitem{Mal06} Malet, A.: Renaissance notions of number and
magnitude.  {\em Historia Mathematica\/} \textbf{33} (2006), no.~1,
63--81.


\bibitem{Mc23} McClenon, R. B.: A Contribution of Leibniz to the
History of Complex Numbers.  {\em American Mathematical Monthly\/}
\textbf{30} (1923), no.~7, 369-374.

\bibitem{Me} M\'eray, H. C. R.: Remarques sur la nature des
quantit\'es d\'efinies par la condition de servir de limites \`a des
variables donn\'ees, \emph{Revue des soci\'eti\'es savantes des
d\'epartments, Section sciences math\'ematiques, physiques et
naturelles, 4th ser.}, \textbf{10} (1869), 280--289.

\bibitem{Nae} Naets, J.: How to define a number? A general
epistemological account of Simon Stevin's art of defining.
\emph{Topoi} \textbf{29} (2010), no.~1, 77--86.


\bibitem{Ro66} Robinson, A.: Non-standard analysis.  North-Holland
Publishing Co., Amsterdam 1966.

\bibitem{She87} Sherry, D.: The wake of Berkeley's Analyst: {\em rigor
mathematicae\/}? {\em Stud. Hist. Philos. Sci.\/} \textbf{18} (1987),
no.~4, 455--480.


\bibitem{St85} Stevin, S. (1585) L'Arithmetique.  In Girard, A. (Ed.),
Les Oeuvres Mathematiques de Simon Stevin (Leyde, 1634), part I,
p.~1--101.


\bibitem{Ste} Stevin, S.: The principal works of Simon
Stevin. Vols. IIA, IIB: Mathematics. Edited by D. J. Struik
C. V. Swets \& Zeitlinger, Amsterdam 1958.  Vol. IIA: v+pp. 1--455
(1~plate).  Vol. IIB: 1958 iv+pp. 459--976.

\bibitem{Tan3} Tannery, P., Henry, C.: {\em Oeuvres de Fermat,
Vol. 3\/} Gauthier-Villars, 1891.



\end{thebibliography}
\end{document}